\documentclass[12pt]{article}
\usepackage{amsmath,amssymb,amsthm,amscd,graphicx}
\usepackage{latexsym}
\usepackage{enumerate}
\usepackage{amsbsy, amsfonts, textcomp}
\usepackage[flushmargin]{footmisc}
\newtheorem{theorem}{Theorem}[]
\newtheorem{proposition}[theorem]{Proposition}

\newtheorem{corollary}[theorem]{Corollary}

\theoremstyle{definition}
\newtheorem{definition}[theorem]{Definition}

\theoremstyle{remark}

\newcommand{\cL}{\mathcal{L}}

\newcommand{\oC}{\overline{C}}
\newcommand{\oK}{\overline{K}}
\newcommand{\oL}{\overline{L}}
\newcommand{\oE}{\overline{E}}
\newcommand{\oF}{\overline{F}}
\def \GL{\operatorname{GL}}
\def \Gal{\operatorname{Gal}}

\def \DGal{\operatorname{DGal}}

\begin{document}

\Large \begin{center} {\bf Galois correspondence theorem}

{\bf for Picard-Vessiot extensions}

\end{center}

\large \centerline{Teresa Crespo, Zbigniew Hajto, El\.zbieta Sowa-Adamus} \normalsize
\let\thefootnote\relax\footnotetext{T. Crespo and Z. Hajto acknowledge support of grant MTM2012-33830, Spanish Science Ministry. E. Sowa-Adamus acknowledges support of the Polish Ministry of Science and Higher
Education. }%%%
\begin{abstract}
In this paper, we generalize the definition of the differential
Galois group and the Galois correspondence theorem established
previously for Picard-Vessiot extensions of real differential
fields with real closed field of constants to any Picard-Vessiot
extension.
\end{abstract}

\section{Introduction}

For a homogeneous linear differential equation $\mathcal{L}(Y)=0$
defined over a differential field $K$ with field of constants $C$,
a Picard-Vessiot extension is a differential field $L$,
differentially generated over $K$ by a fundamental system of
solutions of $\mathcal{L}(Y)=0$ and with constant field equal to
$C$. A classical result states that the Picard-Vessiot extension
exists and is unique up to $K$-differential isomorphism in the
case $C$ algebraically closed  (see \cite{kol2}). Recently, an
existence and uniqueness result for Picard-Vessiot extensions has
been established in the case when the differential field $K$ is a
formally real field (resp. a formally $p$-adic field) with real
closed (resp. $p$-adically closed of the same rank than $K$) field
of constants $C$ (\cite{chp}). In \cite{chs} we presented  a
Galois correspondence theorem for Picard-Vessiot extensions of
formally real differential fields with real closed field of
constants. In this paper we establish a Galois correspondence
theorem for general Picard-Vessiot extensions, valid in particular
for Picard-Vessiot extensions of formally $p$-adic differential
fields with $p$-adically closed field of constants.

Kolchin introduced the concept of strongly normal differential
field extension and obtained a satisfactory Galois correspondence
theorem for this class of extensions without assuming the field of
constants of the differential base field to be algebraically
closed (see \cite{kol} Chapter VI). Note that, for a strongly
normal extension $L|K$, in the case when the constant field of $K$
is not algebraically closed, the differential Galois group is no
longer the group $DAut_K L$ of $K$-differential automorphisms of
$L$, rather one has to consider as well $K$-differential morphisms
of $L$ in larger differential fields. It is worth noting that a
Picard-Vessiot extension is always strongly normal. In the case of
Picard-Vessiot extensions, we can adopt a definition of the
differential Galois group inspired by Kolchin's but simpler than
his one. We obtain then a Galois correspondence theorem which
classifies intermediate differential fields of a Picard-Vessiot
extension in terms of its differential Galois group.

In this paper, we shall deal with fields of
characteristic 0. For the sake of simplicity in the exposition we consider ordinary differential fields.

\section{Main result}

We recall now the precise definition of Picard-Vessiot extension.

\begin{definition}\label{PV} Given a homogeneous linear differential equation

$$\mathcal{L}(Y):=Y^{(n)}+a_{n-1}Y^{(n-1)}+ \ldots + a_{1}Y'+a_{0}Y=0$$

\noindent of order $n$ over a differential field $K$, a
differential extension $L|K$ is a \emph{Picard-Vessiot extension}
for $\cL$ if
\begin{enumerate}
\item $L=K\langle \eta_1,\dots,\eta_n\rangle$, where $\eta_1,\dots,\eta_n$ is a fundamental
set of solutions of $\cL(Y)=0$  in $L$.
\item Every constant of $L$ lies in $K$, i.e. $C_K=C_L$.
\end{enumerate}
\end{definition}

As mentioned in the introduction, a Picard-Vessiot extension is
strongly normal. Hence, the fundamental theorem established by
Kolchin in \cite{kol} chapter VI applies to Picard-Vessiot
extensions. However, for a strongly normal extension $L|K$,
Kolchin defines the differential Galois group $DGal(L|K)$ by means
of differential $K$-isomorphisms of $L$ in the differential
universal extension $U$ of $L$. He obtains then that $DGal(L|K)$
has the structure of an algebraic group defined over the field of
constants $C_U$ of $U$. Afterwards, by using the notion of specialization in Weil's
algebraic geometry, he proves that there exists an algebraic
group $G$ defined over the field of constants $C$ of $K$ such that
$G(C_U)=DGal(L|K)$ . In this section, we give a more direct
definition of the differential Galois group of a Picard-Vessiot
extension, we endow it with a linear algebraic group structure
over $C$ and establish a Galois correspondence theorem in our
setting.

\subsection{Differential Galois group}

Let $K$ be a differential field with field of constants $C$, let
$\mathcal{K}$ be a differential closure of $K$. Let $\oC$ denote
an algebraic closure of $C$ contained in $\mathcal{K}$ and let $\{
\alpha_i \}_{i \in I}$ be a $C$-basis of $\oC$. In the sequel,
differential field extensions of $K$ will be assumed to be
subfields of $\mathcal{K}$. For any such extension $F$ (included
$F=K$), we shall denote by $\oF$ the composition field of $F$ and
$\oC$ inside $\mathcal{K}$.  If the field of constants of $F$ is
equal to $C$, the extensions $\oC|C$ and $F|C$ are linearly
disjoint and then $\{ \alpha_i \}_{i \in I}$ is an $F$-basis of
$\oF$. For a Picard-Vessiot extension $L|K$, we shall consider the
set $DHom_K(L,\oL)$ of $K$-differential morphisms from $L$ into
$\oL$. We shall see that we can define a group structure on this
set and we shall take it as the differential Galois group
$\DGal(L|K)$ of the Picard-Vessiot extension $L|K$. We shall prove
that it is a $C$-defined (Zariski) closed subgroup of some
$\oC$-linear algebraic group.

We observe that we can define mutually inverse bijections

$$\begin{array}{ccc} DHom_K(L,\oL) & \rightarrow & DAut_{\oK}\oL \\ \sigma & \mapsto & \widehat{\sigma}  \end{array},
\quad \begin{array}{ccc}  DAut_{\oK}\oL & \rightarrow &
DHom_K(L,\oL)
\\ \tau & \mapsto & \tau_{|L} \end{array},$$

\noindent where $\widehat{\sigma}$ is the extension of $\sigma$ to
$\oL$. For an element $\sum \lambda_i \alpha_i$ in $\oL$, where
$\lambda_i \in L$, we define $\widehat{\sigma}(\sum \lambda_i
\alpha_i)=\sum \sigma(\lambda_i) \alpha_i$. We may then transfer
the group structure from $DAut_{\oK}\oL$ to $DHom_K(L,\oL)$. Let
us note that $DAut_{\oK}\oL$ is the differential Galois group of
the Picard-Vessiot extension $\oL|\oK$.

Let now $\eta_1,\dots,\eta_n$ be $C$-linearly independent elements
in $L$ such that $L=K\langle \eta_1,\dots,\eta_n\rangle$ and
$\sigma \in DHom_K(L,\oL)$. We have then
$\sigma(\eta_j)=\sum_{i=1}^n c_{ij} \eta_i,$  $1\leq j \leq n$,
with $c_{ij} \in \overline{C}$. We may then associate to $\sigma$
the matrix $(c_{ij})$ in $\GL(n,\overline{C})$. The proofs of
Propositions 16 and 17 and Corollary 18 in \cite{chs} remain valid
in our present setting. We obtain then the following results.

\begin{proposition} Let $K$ be a  differential field with  field of
constants $C$, $L=K\langle \eta_1,\dots,\eta_n\rangle$ a
Picard-Vessiot extension of $K$, where $ \eta_1,\dots,\eta_n$ are
$C$-linearly independent. There exists a set $S$ of polynomials
$P(X_{ij}), 1\leq i,j \leq n$, with coefficients in $C$ such that
\begin{enumerate}[1)]
\item If $\sigma \in DHom_K(L,\oL)$ and
$\sigma(\eta_j)=\sum_{i=1}^n c_{ij}\eta_i$, then $P(c_{ij})=0, \forall P
\in S$.
\item Given a matrix $(c_{ij}) \in \GL(n,\overline{C})$ with $P(c_{ij})=0, \forall P
\in S$, there exists a differential $K$-morphism $\sigma$ from $L$
to $\oL$ such that $\sigma(\eta_j)=\sum_{i=1}^n c_{ij}\eta_i$.
\end{enumerate}
\end{proposition}

The preceding proposition gives that $\DGal(L|K)$ is a $C$-defined
closed subgroup of $\GL(n,\overline{C})$.

\begin{proposition}
Let $K$ be a differential field with field of constants $C$, $L|K$
a  Picard-Vessiot extension. For $a \in L\setminus K$, there
exists a $K$-differential morphism $\sigma:L \rightarrow \oL$ such
that $\sigma(a) \neq a$.
\end{proposition}

\vspace{0.5cm}
For a subset $S$ of $\DGal(L|K)$, we set $L^S:=\{ a \in L : \sigma(a)=a, \, \forall \sigma \in S \}$.

\begin{corollary}\label{cor}
Let $K$ be a  differential field with field of constants $C$,
$L|K$ a  Picard-Vessiot extension. We have  $L^{\DGal(L|K)}=K$.
\end{corollary}

\subsection{Fundamental theorem}\label{tf}

Let $K$ be a differential field with field of constants $C$ and
$L|K$ a  Picard-Vessiot extension. For a closed subgroup $H$ of
$\DGal(L|K)$, $L^H$ is a differential subfield of $L$ containing
$K$. If $E$ is an intermediate differential field, i.e. $K\subset
E \subset L$, then $L|E$ is a  Picard-Vessiot extension and
$\DGal(L|E)$ is a $C$-defined closed subgroup of $\DGal(L|K)$.

\begin{theorem} Let $L|K$ be a Picard-Vessiot extension, $\DGal(L|K)$ its differential Galois
group.
\begin{enumerate}[1.]
\item
The correspondences

$$H\mapsto L^H \quad , \quad E\mapsto \DGal(L|E)$$

\noindent define inclusion inverting mutually inverse bijective
maps between the set of $C$-defined closed subgroups $H$ of $\DGal(L|K)$
and the set of differential fields $E$ with $K\subset E \subset
L$.

\item The intermediate differential field $E$ is a Picard-Vessiot extension of $K$ if and only if the subgroup
$\DGal(L|E)$ is normal in $\DGal(L|K)$. In this case, the restriction morphism

    $$\begin{array}{ccc} \DGal(L|K)& \rightarrow & \DGal(E|K) \\ \sigma & \mapsto & \sigma_{|E} \end{array} $$

\noindent induces an isomorphism $$\DGal(L|K)/\DGal(L|E) \simeq \DGal(E|K).$$
\end{enumerate}
\end{theorem}

\noindent {\it Proof.} 1.It is clear that both maps invert
inclusion. If $E$ is an intermediate differential field of $L|K$,
we have $L^{\DGal(L|E)}=E$, taking into account that $L|E$ is
Picard-Vessiot and corollary \ref{cor}. For $H$ a $C$-defined
closed subgroup of $\DGal(L|K)$, the equality $H=\DGal(L|L^{H})$
follows from the correspondent equality in Picard-Vessiot theory
for differential fields with algebraically closed field of
constants (see e.g. \cite{haj} theorem 6.3.8).

\noindent 2.If $E$ is a Picard-Vessiot extension of $K$, then
$\oE$ is a Picard-Vessiot extension of $\oK$ and so $\DGal(L|E)$
is normal in $\DGal(L|K)$. Reciprocally, if $\DGal(L|E)$ is normal
in $\DGal(L|K)$, then the subfield of $\oL$ fixed by $\DGal(L|E)$
is a Picard-Vessiot extension of $\oK$. Now, this field is $\oE$.
So, $\oE$ is differentially generated over $\oK$ by a
$\overline{C}$-vector space $V$ of finite dimension. Let
$\{v_1,\dots,v_n \}$ be a $\oC$-basis of $V$. We may write each
$v_j, 1 \leq j \leq n$ as a linear combination of the elements
$\alpha_i$ with coefficients in $E$. Now, there is a finite number
of $\alpha_i$'s appearing effectively in these linear
combinations. Let $\widetilde{C}$ be a finite Galois extension of
$C$ containing all these $\alpha_i$'s. We have then $v_i \in
\widetilde{E}:=\widetilde{C} \cdot E$ and $\widetilde{E}$ is
differentially generated over $\widetilde{K}:=\widetilde{C} \cdot
K$ by $V$. We may extend the action of $\Gal(\widetilde{C}|C)$ to
$\widetilde{E}$ and consider the translate $c(V)$ of $V$ by $c \in
\Gal(\widetilde{C}|C)$. Let $\widetilde{V}= \oplus_{c \in
\Gal(\widetilde{C}|C)} c(V)$. We have that  $\widetilde{E}$ is
differentially generated over $\widetilde{K}$ by $\widetilde{V}$
and $\widetilde{V}$ is $\Gal(\widetilde{C}|C)$-stable, hence $E$
is differentially generated over $K$ by the $C$-vector space
$\widetilde{V}^{\Gal(\widetilde{C}|C)}=\{y \in V: c(y)=y, \,
\forall c \in \Gal(\widetilde{C}|C)\}$. Hence $E|K$ is a
Picard-Vessiot extension. The last statement of the theorem
follows from the fundamental theorem of Picard-Vessiot theory in
the case of algebraically closed fields of constants (\cite{haj}
theorem 6.3.8). \hfill $\Box$

\vspace{1cm}
\footnotesize
\noindent
Teresa Crespo, Departament d'\`{A}lgebra i Geometria, Universitat de Barcelona, Gran Via de les Corts Catalanes 585, 08007 Barcelona, Spain, teresa.crespo@ub.edu

\vspace{0.2cm}
\noindent
Zbigniew Hajto, Faculty of Mathematics and Computer Science, Jagiellonian University, ul. Prof. S. \L ojasiewicza 6, 30-348 Krak\'ow, Poland, zbigniew.hajto@uj.edu.pl

\vspace{0.2cm}
\noindent
El\.zbieta Sowa-Adamus, Faculty of Applied Mathematics, AGH University of Science and Technology,
al. Mickiewicza 30, 30-059 Krak\'ow, Poland, esowa@agh.edu.pl

\end{document}